\documentclass{nseJournal}

\newcommand{\rmd}{\mathrm{d}}
\usepackage{amsmath,amssymb,amsfonts,amsthm,subcaption,graphicx,url}

\newtheorem{theorem}{Theorem}[]
\newtheorem{corollary}{Corollary}[]

\begin{document}

\title{A Stochastic Calculus Approach to Boltzmann Transport}

\addAuthor{\correspondingAuthor{J.~Darby Smith}}{a}
\correspondingEmail{jsmit16@sandia.gov}
\addAuthor{Rich Lehoucq}{a}
\addAuthor{Brian Franke}{b}

\addAffiliation{a}{Sandia National Laboratories, Center for Computing Research\\P.O. Box 5800, Albuquerque, NM 87123}
\addAffiliation{b}{Sandia National Laboratories, Radiation Effects Theory Department\\ P.O. Box 5800, Albuquerque, New Mexico 87123}

\addKeyword{Boltzmann transport}
\addKeyword{Stochastic differential equations}
\addKeyword{Adjoint transport}

\titlePage

\begin{abstract}
Traditional Monte Carlo methods for particle transport utilize source iteration to express the solution, the flux density, of the transport equation as a Neumann series. Our contribution is to show that the particle paths simulated within source iteration are associated with the adjoint flux density and the adjoint particle paths are associated with the flux density. We make our assertion rigorous through the use of stochastic calculus by representing the particle path used in source iteration as a solution to a stochastic differential equation (SDE). The solution to the adjoint Boltzmann equation is then expressed in terms of the same SDE and the solution to the Boltzmann equation is expressed in terms of the SDE associated with the adjoint particle process. An important consequence is that the particle paths used within source iteration simultaneously provide Monte Carlo approximations of the flux density and adjoint flux density. Monte Carlo simulations are presented to support the simultaneous use of the particle paths.
\end{abstract}

\section{Introduction}

Introduced by Ludwig Boltzmann for the kinetics of particular gasses, rigorous simulation of the Boltzmann equation for particle and radiation transport is critical for power distribution problems, radiation detector design, x-ray simulation, spectral response and analysis, and numerous other radiation applications \cite{dunn2009monte,vaz2009neutron,zhang2022variational}.

While deterministic numerical solvers for certain forms of the Boltzmann transport equation do exist \cite{lewis1984computational,ADAMS20023}, Monte Carlo approaches for numerical approximation remain the dominant approach to solving neutron transport problems \cite{dunn2009monte,vaz2009neutron,peraud2014monte}. There are many variations of Monte Carlo approaches for neutron transport, including those using multilevel schemes \cite{louvin2017adaptive} or employing additional sampling or counting sub-steps like the predictor-corrector method \cite{leppanen2009two}.

Though there are several ways to alter the Monte Carlo approach, the basic method is one of sampling particle paths and ``scoring'' or ``tallying'' them in detector regions \cite{dupree2012monte,lux2018monte}. Particles are typically initialized according to a prescribed source term. When calculated in this mode, particle paths can be reused to calculate new tallies based on different detectors. However, as the method of initialization depends upon the source, the particle paths cannot be reused if the source term changes. Nonetheless, if the detector area is small, many paths simulated from the source may never reach the detector \cite{kersting2020single} so that the adjoint flux density can be of more utility than the flux density. The adjoint simulation is accomplished by initializing particles from the detector and simulating them back to the source, in contrast to the forward simulation \cite{kersting2020single,vitali2018comparison}. By starting at a detector region, with the detector as the adjoint source term, more paths are likely to be useful by reaching the source and the adjoint paths can be reused when the source is modified.  Sometimes, both forward and adjoint information is desired. 
Unfortunately since simulations can have differing scattering rates---sometimes handled through altered rates or through pre-calculated tally weights and likely have differing starting locations---forward transport particle paths cannot be reused to calculate adjoint information, and vice versa, without altering or reweighting the paths. 

Underpinning the use of the flux density or the adjoint flux density embodies the well-known adjoint relationship that the response functional is the inner product of the flux density with the adjoint source or the inner product of the adjoint flux density with the source.
Our contribution is to show that the particle paths simulated within source iteration are associated with the adjoint flux density and the adjoint particle paths are associated with the flux density.
We make our assertion rigorous through the use of stochastic calculus by representing the particle path used in source iteration as a solution to a stochastic differential equation (SDE). The solution to the adjoint Boltzmann equation is then expressed in terms of the same SDE and the solution to the Boltzmann equation is expressed in terms of the SDE associated with the adjoint particle process. An important consequence is that the particle paths used within source iteration simultaneously provide Monte Carlo approximations of the flux density and adjoint flux density in the detector and source regions, respectively. To the best our knowledge, this has not appeared in the literature.

We describe Boltzmann transport in the case of non-fissionable material through an SDE in Section \ref{sec:boltz}, connecting this representation to the adjoint in Section \ref{sec:adjoint_sde}. We follow this with a brief commentary on an analogous result discussing the adjoint process as an SDE in Section \ref{sec:adjoint_theorem}. Then, we provide numerical examples in Section \ref{sec:connecting} and a discussion in Section \ref{sec:discussion}. We provide proofs of our mathematical theorems in the Appendix.


\section{Methods}

\subsection{Boltzmann Transport as an SDE}
\label{sec:boltz}

Traditional transport methods use source iteration to approximate the solution of the Boltzmann equation, i.e., the flux density, in terms of a Neumann series. Each term of the series has an interpretation as a sum over particles that undergo a number of events \cite{dupree2012monte}. In practice, the Neumann series is approximated via Monte Carlo where average tallies are calculated using the history of particles initialized over some region with energy and direction based upon the source term $f$.  A time to interaction is drawn, or, as particles travel with fixed velocity, a distance to interaction is drawn. Tallies for each location the particle passes through on the way to the interaction are increased.  In the case of non-fissionable material, there are two primary types of interactions, scattering and absorption.\footnote{Fission and the branching process do not fit into this framework as presented; however, this is the focus of continuing work.} At the interaction point, a particle either receives a new direction and energy due to scattering or is absorbed. Because absorption and scattering are independent, often only time or distance to scattering are drawn and tallies are all weighted by the probability that the particle would not have been absorbed at that time or location. A particle's journey ends upon exiting the domain $\mathcal{D}$. The preceding description corresponds to a conventional Monte Carlo approach to approximate the solution of Boltzmann equation via source iteration.

We recast the particle trajectories that form the basis for source iteration using stochastic calculus. In order to simplify notation, describe a three dimensional version of Boltzmann transport undergoing absorption and scattering and remark that our presentation is general and can be seen to apply to transport beyond three dimensions. We denote the position of a particle by $X_1(t)$ and its direction and energy by $X_2(t)$ and $X_3(t)$, all three indexed by time. The components of the vector $X(t) = \left(X_1(t),X_2(t),X_3(t)\right)$ are stochastic processes indexed by time but can instead be written in terms of distance traveled with a change of coordinates. The stochastic differential equation (SDE) governing the process $X(t)$  is given by 
\begin{align}
\begin{split}
\rmd X_1(t) & = X_2(t)\, \rmd t,\\
\rmd X_2(t) & = \left(\omega - X_2(t)\right)\rmd N\left(t,\omega, \varepsilon, X(t)\right),\\
\rmd X_3(t) & = \left(\varepsilon - X_3(t)\right)\rmd N\left(t, \omega, \varepsilon, X(t)\right)
\end{split}
\label{eq:adj_sde}
\end{align}
where scattering is given by the Poisson process $N(t,\omega, \varepsilon, X(t))$ with rate $\Sigma_s\left(X(t)\right)$,  and
$\omega$ and $\varepsilon$ are random variables representing post-scatter direction and energy. They are distributed with respect to the conditional probability density $ p\left(\cdot\,\middle|\,X_2(t),X_3(t)\right)$ dependent upon the pre-scatter values. The initial particle position $X(0)$ is distributed according to the source $f$, which we assume to be a probability density and the SDE \eqref{eq:adj_sde} holds until the particle exits the domain $\mathcal{D}$. The exit-time, a random variable and function of the starting location $(x,\Omega,E)$, is given by 
\begin{equation}
\tau_{x,\Omega,E}  = \inf \left\{t>0\,\middle|\,X(0)=\left(x,\Omega,E\right),X(t)\not\in\mathcal{D}\right\}.
\label{eq:stopping_time}
\end{equation}
In absorbing mediums, particles may also be randomly absorbed with exponentially distributed absorption times with rate $\Sigma_a\left(X(t)\right)$. In this framework, absorption is not explicitly accounted for in the SDE and is instead handled within the expectation expression for the Boltzmann flux density. The SDE representation given by \eqref{eq:adj_sde} is understood to be a compact version of 
\begin{align}
\begin{split}
X_1(t) & = X_1(0) + \int_0^tX_2(s)\,\rmd s,\\
X_2(t) & = X_2(0) + \int_0^{t} (\omega - X_2(s) ) \,\rmd N(s,\omega, \varepsilon, X(s))\\
X_3(t) & = X_3(0) + \int_0^{t}(\varepsilon - X_3(s)) \, \rmd N (s, \omega, \varepsilon, X(s))
\end{split}
\label{eq:adj_sde-full}
\end{align}
where the quantity $\rmd N$ is the counting process associated with the Poisson process $N$ to indicate that an event occurs. The counting process $\rmd N$ may be identified with a Dirac delta function, which is zero except at the event times induced by the  Poisson process $N$ so that the integrals reduce to sums. The equation for $X_1$ explains that the particle transports until the particle exits or scatters. 
Following a scattering event, the process $X_2$ is updated by a random change in direction $\omega-X_2(t)$ and in an analogous fashion for $X_3$. Hence, when a Poisson event occurs, the direction of the process is replaced by the random value $\omega$ and the energy by $\varepsilon$.

The inner product of the flux density $\Phi$ and a detector $g$, i.e.,  $\langle g, \Phi\rangle$  can now be expressed as the expectation
\begin{equation}
\langle g,\Phi\rangle = \mathbb{E}_{f}\left[\mathbb{E}\left[\int_0^{\tau_{x,\Omega,E}}g\left(X(t)\right)\exp\left(-\int_0^t\Sigma_a\left(X(u)\right)\rmd u\right)\rmd t\,\middle|\,X(0)=\left(x,\Omega,E\right)\right]\right].
\label{eq:trad_boltz_sde}
\end{equation}
The integral accumulates tallies based upon the detector  $g$ until the exit-time and are attenuated by the probability that a particle has not been absorbed by the time it enters the domain of $g$. The inner (conditional) expectation is over particles starting at the point $\left(x,\Omega,E\right)$ and the outer expectation is taken so that starting points are chosen and averaged with respect to the source $f$.

\subsection{Relating the Adjoint to the SDE}
\label{sec:adjoint_sde}

Let $\Psi$ denote the adjoint flux density so that the adjoint relationship explains that the response functional is given by the equality
\begin{equation}
\langle \Psi, f\rangle = \langle g, \Phi \rangle
\label{eq:adjoint_property}
\end{equation}
where the detector $g$ is the source term for the adjoint transport equation.
The left-most inner product is an integral of $\Psi$ against a density, and is therefore an expectation so that inserting \eqref{eq:trad_boltz_sde} into \eqref{eq:adjoint_property} results in 
\begin{equation}
\mathbb{E}_{f}[\Psi] = \langle \Psi, f\rangle =  \mathbb{E}_{f}\left[\mathbb{E}\left[\int_0^{\tau_{x,\Omega,E}}g\left(X(t)\right)\exp\left(-\int_0^t\Sigma_a\left(X(u)\right)\rmd u\right)\rmd t\,\middle|\,X(0)=\left(x,\Omega,E\right)\right]\right].
\label{eq:adj_flux_density}
\end{equation}

We now state a theorem that asserts conditions under which we can identify the adjoint flux density with the inner expectation of \eqref{eq:adj_flux_density}. That is, the adjoint flux density is a sum of conditional expectations, one of which is a function $g$ of the solution $X$ of the SDE \eqref{eq:adj_sde} and the other involves a prescribed boundary condition for the adjoint flux density. 
The proof of the theorem employs the stochastic calculus and is given in the Appendix.

\begin{theorem}\label{thm:adj}
Let $X(t)$ be given by \eqref{eq:adj_sde} and $\tau_{x,\Omega,E}$ be given by \eqref{eq:stopping_time}. Let $\Sigma_t=\Sigma_a+\Sigma_s$.  Suppose there exists a classical solution $\Psi\colon\mathcal{D}\to \mathbb{R}$ where $\mathcal{D}$ is a compact set for the adjoint Boltzmann equation
\begin{align}
\begin{split}
-\Omega \frac{\partial}{\partial x}\Psi\left(x,\Omega,E\right) &+ \Sigma_t\left(x,\Omega,E\right)\Psi\left(x,\Omega,E\right) - g\left(x,\Omega,E\right)\\
&= \int\Psi\left(x,\Omega',E'\right)\Sigma_s\left(x,\Omega,E\right)p\left(\Omega',E'\,\middle|\,\Omega,E\right)\rmd E'\rmd\Omega',\\
\Psi\left(x,\Omega,E\right) &=B\left(x,\Omega,E\right),\qquad \left(x,\Omega,E\right)\in\partial\mathcal{D}.\\
\end{split}
\label{eq:adj_bvp}
\end{align}
Then, when $\Sigma_a$, $\Sigma_s$, and $g$ are continuous almost everywhere and bounded, and if $\mathbb{E}\left[\tau_{x,\Omega,E}\right]<\infty$ for all $(x,\Omega,E)\in\mathcal{D}$, then the adjoint flux density $\Psi$ is
\begin{align}
\begin{split}
\Psi\left(x,\Omega,E\right) = &\mathbb{E}\left[B\left(X\left(\tau_{x,\Omega,E}\right)\right)\exp\left(-\int_0^{\tau_{x,\Omega,E}}\Sigma_a\left(X(t)\right)\,\rmd t\right)\,\middle|\, X(0)=\left(x,\Omega,E\right)\right]\\
&+ \mathbb{E}\left[\int_0^{\tau_{x,\Omega,E}} g\left(X(t)\right)\exp\left(-\int_0^t\Sigma_a\left(X(u)\right)\,\rmd u\right)\,\rmd t\,\middle|\, X(0)=\left(x,\Omega,E\right)\right].
\end{split}
\label{eq:adj_sol_general}
\end{align}
\end{theorem}

The significant conclusion is that we can simultaneously use the particle trajectories, i.e., samples of the process $X$ generated during source iteration to approximate the
\begin{enumerate}
\item flux density $\Phi$ over the source region, and 
\item adjoint flux density $\Psi$ over the detector region.
\end{enumerate}
The use of the particle paths for the first approximation is well-understood but the second use of the paths in this manner does not appear in the literature. Our use of the stochastic calculus establishes this remarkable dual use of the particle paths. 
The practical consequence is significant because standard practice involving the adjoint flux density initializes an adjoint particle path with respect to a detector $g$ so preventing the reuse of those trajectories for distinct detectors.

In the common case of vacuum boundary conditions, Theorem \ref{thm:adj} explicitly states that the adjoint flux density is given as the expectation of particles scored in the detector $g$, weighted by an absorption probability:
\begin{equation}
\Psi\left(x,\Omega,E\right)=\mathbb{E}\left[\int_0^{\tau_{x,\Omega,E}}g\left(X(t)\right)\exp\left(-\int_0^t\Sigma_a\left(X(u)\right)\rmd u\right)\rmd t\,\middle|\,X(0)=\left(x,\Omega,E\right)\right].
\label{eq:adj_solution}
\end{equation}


\subsection{Analogous Results Utilizing the Adjoint Process}
\label{sec:adjoint_theorem}

The goal of this section is to state a Corollary to Theorem \ref{thm:adj} that explains that the flux density is a conditional expectation over the particle trajectories used within adjoint source iteration. 

We first describe the adjoint particle process as an SDE where the particles are initialized with respect to adjoint source (detector) $g$, which we assume to be a probability density. The adjoint particles move in a manner that is reversed from the traditional process $X(t)$. Scattering rates and distributions for the adjoint process are calculated from the original scattering rate $\Sigma_s$ and scattering distribution $p$. In line with our established notation, we write the adjoint scattering cross section as
\begin{equation}
S \left(x,\Omega,E\right) = \int \Sigma_s\left(x,\Omega',E'\right)p\left(\Omega,E\,\middle|\,\Omega',E'\right)\,\rmd\Omega'\rmd E.
\label{eq:back_scatter}
\end{equation}
The adjoint scattering distribution $q$ is defined so that balance is achieved with respect to both rates:
\begin{equation}
S\left(x,\Omega,E\right)q\left(\Omega',E'\,\middle|\,\Omega,E\right) = \Sigma_s\left(x,\Omega',E'\right)p\left(\Omega,E\,\middle|\,\Omega',E'\right).
\label{eq:q_def}
\end{equation}

We will now define an adjoint particle process $Y(t) = \left(Y_1(t),Y_2(t),Y_3(t)\right)$ in one spatial dimension analogous to \eqref{eq:adj_sde}. The SDE is
\begin{align}
\begin{split}
\rmd Y_1(t) & = -Y_2(t) \, \rmd t,\\
\rmd Y_2(t) & = (\zeta- Y_2(t)) \, \rmd M\left(t,\zeta, \rho, Y(t)\right),\\
\rmd Y_3(t) & = (\rho - Y_3(t)) \, \rmd M\left(t, \zeta, \rho, Y(t)\right),
\end{split}
\label{eq:boltz_sde}
\end{align}
and the associated stopping time is 
\begin{equation}
\sigma_{x,\Omega,E}  = \inf \left\{t>0\,\middle|\,Y(0)=\left(x,\Omega,E\right),Y(t)\not\in\mathcal{D}\right\}.
\label{eq:adj_particle_stopping_time}
\end{equation}
Here, the position of an adjoint particle is $Y_1(t)$ and its direction and energy are $Y_2(t)$ and $Y_3(t)$. The components of the vector $Y(t) = \left(Y_1(t),Y_2(t),Y_3(t)\right)$ are stochastic processes indexed by time but can instead be written in terms of distance traveled with a change of coordinates. 
The particle scatters according to a Poisson process $M\left(t,\zeta,\rho,Y(t)\right)$ with rate $S\left(Y(t)\right)$ until the exit-time $\sigma_{x,\Omega,E}$ from the domain $\mathcal{D}$. When scattering occurs, a new direction $\zeta$ and a new energy $\rho$ are selected given the pre-scatter state from the conditional distribution $q$ defined in \eqref{eq:q_def}.  Once a scattering event occurs, both the direction $Y_2$ and the energy $Y_3$ are updated accordingly. 

The SDE representation given by \eqref{eq:boltz_sde} is understood to be a compact version of 
\begin{align}
\begin{split}
Y_1(t) & = Y_1(0) -  \int_0^tY_2(s)\,\rmd s,\\
Y_2(t) & = Y_2(0) + \int_0^{t} (\zeta- Y_2(s)) \,\rmd M(s,\zeta, \rho, Y(s))\\
Y_3(t) & = Y_3(0) + \int_0^{t} (\rho - Y_3(s))\, \rmd M (t, \zeta, \rho, Y(s))
\end{split}
\label{eq:boltz_sde-full}
\end{align}
where the quantity $\rmd M$ is the counting process associated with the Poisson process $M$ to indicate that an event occurs. The counting process $\rmd M$ may be identified with a Dirac delta function, which is zero except at the event times induced by the  Poisson process $M$ so that the integrals reduce to sums.

The inner product of the adjoint flux density $\Psi$ and the source $f$, i.e.,  $\langle f, \Psi\rangle$  can now be expressed as the expectation
\begin{equation}
\langle f,\Psi\rangle = \mathbb{E}_{g}\left[\mathbb{E}\left[\int_0^{\tau_{x,\Omega,E}}f\left(Y(t)\right)\exp\left(-\int_0^t(\Sigma_t-S)(Y(u))\rmd u\right)\rmd t\,\middle|\,Y(0)=\left(x,\Omega,E\right)\right]\right]
\label{eq:ave-adj_boltz_sde}
\end{equation}
where $\Sigma_t = \Sigma_a+\Sigma_s$. 
The integral accumulates tallies based upon the source $f$ until the exit-time and are attenuated by the probability that a particle has not been absorbed by the time it enters the domain of $f$. The inner (conditional) expectation is over particles starting at the point $\left(x,\Omega,E\right)$ and the outer expectation is taken so that starting points are chosen and averaged with respect to the detector $g$.

The following corollary of Theorem \ref{thm:adj}  asserts conditions under which we can identify the flux density as a sum of conditional expectations, one of which is a function $f$ of the solution  $Y$ of the SDE \eqref{eq:boltz_sde} and the other involves a prescribed boundary condition for the flux density. 
The proof of this is provided in the Appendix.

\begin{corollary}\label{thm:boltz}
Let $Y(t)$ and $\sigma_{x,\Omega,E}$ be given by \eqref{eq:boltz_sde}, and $\Sigma_t = \Sigma_a+\Sigma_s$. Suppose that there exists a classical solution $\Phi\colon\mathcal{D}\to\mathbb{R}$ where $\mathcal{D}$ is compact that solves the Boltzmann equation
\begin{align}
\begin{split}
\frac{\partial}{\partial x}\Omega \Phi\left(x,\Omega,E\right)&+\Sigma_t\left(x,\Omega,E\right)\Phi(x,\Omega,E) - f\left(x,\Omega,E\right)\\
& = \int\Phi\left(x,\Omega',E'\right)\Sigma_s\left(x,\Omega',E'\right)p\left(\Omega, E\,\middle|\,\Omega',E'\right)\,\rmd E'\rmd\Omega',\\
\Phi\left(x,\Omega,E\right) & = H\left(x,\Omega,E\right),\qquad\left(x,\Omega,E\right)\in \partial \mathcal{D}.
\end{split}
\label{eq:boltz_problem}
\end{align}
Then, when $\Sigma_a$, $\Sigma_s$, and $f$ are continuous almost everywhere and bounded, and if $\mathbb{E}\left[\sigma_{x,\Omega,E}\right]<\infty$ for all $\left(x,\Omega,E\right)\in\mathcal{D}$, then the flux density $\Phi$  is
\begin{align}\small
\begin{split}
\Phi\left(x,\Omega,E\right) = \mathbb{E}&\left[H\left(Y\left(\sigma_{x,\Omega,E}\right)\right)\exp\left(-\int_0^{\sigma_{x,\Omega,E}}\left(\Sigma_t-S\right)(Y(t))\,\rmd t\right)\,\middle|\,Y\left(0\right)=\left(x,\Omega,E\right)\right]\\
&+\mathbb{E}\left[\int_0^{\sigma_{x,\Omega,E}}f\left(Y\left(t\right)\right)\exp\left(-\int_0^t\left(\Sigma_t-S\right)(Y(u))\,\rmd u\right)\rmd t\,\middle|\,Y(0)=\left(x,\Omega,E\right)\right].
\end{split}
\label{eq:boltz_sde_solution}
\end{align}\normalsize
\end{corollary}

The significant conclusion is that we can simultaneously use the adjoint particle trajectories, i.e., samples of the process $Y$ to approximate the
\begin{enumerate}
\item adjoint flux density $\Psi$ over the detector region, and 
\item flux density $\Phi$ over the source region.
\end{enumerate}
The use of the particle paths for first approximation is well-understood but the second use of the paths in this manner does not appear in the literature. Our use of the stochastic calculus establishes this remarkable dual use of the particle paths. 
The practical consequence is significant because standard practice involving the flux density initializes the particle process $X$ with respect to a source $f$ so preventing the reuse of those trajectories for distinct source.


\section{Results}
\label{sec:connecting}

The immediate impact of Theorem \ref{thm:adj} is that trajectories simulated for Boltzmann transport in non-fissionable material can be reused to obtain the adjoint quantity.  Analogously, Corollary \ref{thm:boltz} allows adjoint trajectories simulated for the adjoint quantity to be reused to obtain the Boltzmann transport solution.

We showcase these connections in two ways.  First, in Section \ref{sec:simple_prob}, we demonstrate reuse of trajectories for a simple one-dimensional problem with monotone energy. This simplified problem allows one to immediately compare quantities obtained through reuse of trajectories to their traditionally simulated counterparts.  Then, in Section \ref{sec:its} we demonstrate the reuse of trajectories in the Integrated TIGER Series (ITS) code \cite{franke2023} via a full-dimensional problem (3 spatial dimensions, 2 angular, 1 energy) utilizing both energy and angle in a photon transport problem using realistic physics.

\subsection{Simplified Problem}
\label{sec:simple_prob}

We now demonstrate the reuse of samples in the simplified setting of 2D fluence problem (1 spatial, 1 angle/direction) with no energy dependence. We use the parameterization in Table \ref{tab:params_sim} and vacuum boundary conditions on $\mathcal{D} = [-1,1]\times[-1,1]$.
\begin{table}[h]\centering
\caption{Table of variable values for numerical example. Note, $\chi$ is the indicator function for the given rectangles and $\mathcal{A}$ represents the area of the specified rectangle.}
\begin{tabular}{l|c|c}\hline
\textbf{Variable}                     & \textbf{Value}    & \textbf{Units} \\\hline\hline
$\Sigma_a\left(\Omega\right)$ & $5.00$               &  $\text{cm}^{-1}$\\
$\Sigma_s\left(\Omega\right)$ & $2.50$               &   $\text{cm}^{-1}$\\
$f\left(x,\Omega\right)$          & $\mathcal{A}\left(R_1\right)\chi_{R_1}\left(x,\Omega\right)$ & $\text{cm}^{-3}\text{steradian}^{-1}$  \\
$R_1$                                   & $[0.29,0.69]\times[-1,1]$ & Rectangle in $(x,\Omega)$ space \\
$g\left(x,\Omega\right)$         & $\mathcal{A}\left(R_2\right)\chi_{R_2}\left(x,\Omega\right)$ & $\text{cm}^{-3}\text{steradian}^{-1}$ \\
$R_2$                                   & $[-0.22,-0.06]\times[-1,1]$ & Rectangle in $(x,\Omega)$ space \\
$p\left(\Omega\,\middle|\Omega'\right)$ & Uniform, no dependence on $\Omega'$. & Probability density \\ \hline
\end{tabular}
\label{tab:params_sim}
\end{table}

\begin{figure}[htp]
\centering
\begin{subfigure}{.49\textwidth}
\includegraphics[width=\textwidth]{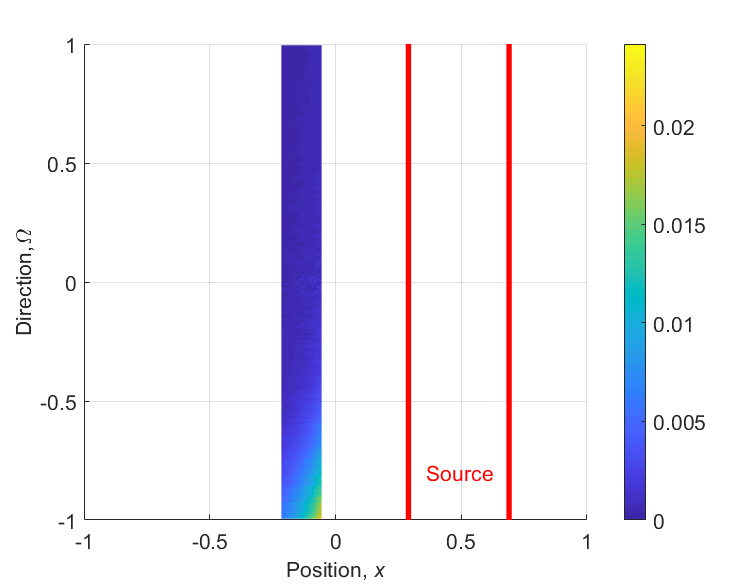}
\caption{Traditional Boltzmann}
\label{fig:boltz_trad}
\end{subfigure}
\begin{subfigure}{.49\textwidth}
\includegraphics[width=\textwidth]{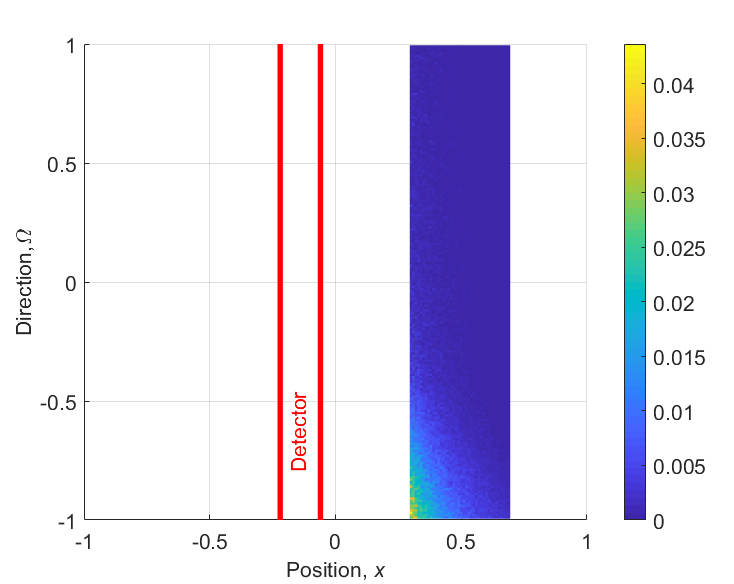}
\caption{Adjoint, Reused Trajectories}
\label{fig:adj_boltz_reuse}
\end{subfigure}
\begin{subfigure}{.49\textwidth}
\includegraphics[width=\textwidth]{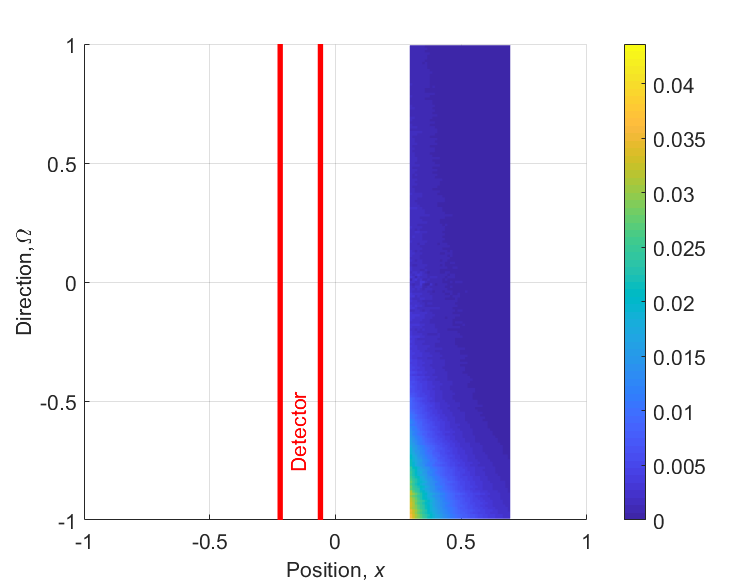}
\caption{Traditional Adjoint}
\label{fig:adj_trad}
\end{subfigure}
\begin{subfigure}{.49\textwidth}
\includegraphics[width=\textwidth]{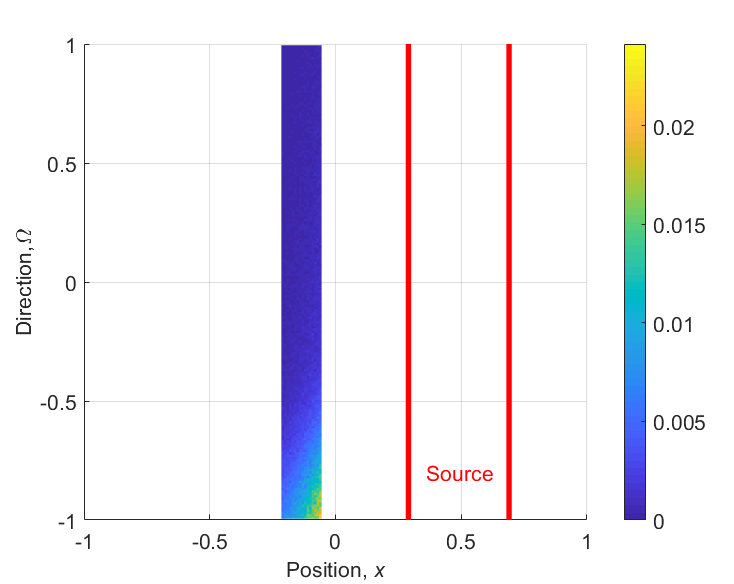}
\caption{Boltzmann, Reused Trajectories}
\label{fig:boltz_reuse}
\end{subfigure}
\caption{\textbf{(a)} Traditional Boltzmann fluence simulation using forward particles. Particles travel from source (marked in red) to detector (area of solution). \textbf{(b)} Adjoint fluence simulation generated by reusing trajectories generated in \textbf{(a)}. Particles travel from the source (area of solution) and are scored in the detector (marked in red). \textbf{(c)} Traditional adjoint fluence simulation using adjoint particles. Particles move from detector (marked in red) to source (area of solution). \textbf{(d)} Boltzmann fluence generated through the reuse of adjoint particles simulated in \textbf{(c)}.  The reuse of particles in \textbf{(b)} generates a good simulation, with the norm of difference between \textbf{(b)} and \textbf{(c)} equal to $0.0443$.  Similarly, the reuse of particles in \textbf{(d)} produces a simulation close to that of \textbf{(a)}, with the norm of the difference equal to $0.0162$.}
\label{fig:numerical_ex}
\end{figure}

We first completed a traditional Boltzmann simulation by initializing a population of $500,200$ particles so that their density approximates the source term $f$. Since our selected source term is flat, we allocated an equal number of particles across a mesh with spatial discretization $\Delta s=0.01$ and direction/angular discretization $\Delta a = 0.01$.  Particles were simulated in discrete time with step size $\Delta t=0.01$. Once initialized, a time to scattering is selected as an exponential random variable with parameter $\Sigma_s^{-1}$. This time is fit to the discretization by selecting the first time step that exceeds this time to event. The particle increments its location by $v\Omega\Delta t$ until this event time step where a new direction $\Omega$ is selected. For each completed time step, a tally is recorded if the particle is within a mesh point of the detector rectangle $R_2$. This tally is weighted by the probability that the particle would not have been absorbed by this time step. This process repeats until the particle leaves the domain $\mathcal{D}$, where we stop tracking it.  The solution calculated over the sensor domain is given in Fig.~\ref{fig:boltz_trad}.

We reused these exact trajectories to calculate the adjoint quantity in the support of $f$. From \eqref{eq:adj_solution}, we can simulate the adjoint quantity at any point $(x,\Omega)$ through the use of particle trajectories initialized at that point. Reusing the trajectories we have already generated, this limits calculation to those mesh points in the support of $f$.  In a discretized fashion, we examine trajectories beginning within a mesh point. Tracing those trajectories, we keep a score, adding the constant value $\mathcal{A}\left(R_2\right)e^{-\Sigma_a j\Delta t}$ to our score whenever the trajectory is in the rectangle $R_2$ on the $j^{\text{th}}$ time step. Once all trajectories for all particles starting within the mesh point have been traced and the score calculated, we divide that score by the total number of particles starting within the mesh point (in the case of this simulation, this is always 61). This division takes the expected value over all particles initialized and gives us the adjoint value for that mesh point.  We plot our solution obtained from the reuse of trajectories in Fig.~\ref{fig:adj_boltz_reuse}.

Next, we performed a traditional adjoint simulation by simulating adjoint particles initialized so that their density approximates the adjoint source/sensor $g$. Aiming to use a similar number of total particles as the traditional Boltzmann transport simulation, we used $499,800$ total particles yielding $147$ particles per location in the initial discretization of $R_2$. The results from our adjoint simulation are in Fig.~\ref{fig:adj_trad}.  Comparing the adjoint results, we immediately see visual agreement though the reused trajectories plot does seem noisier. This is to be expected -- the reuse of trajectories calculated the adjoint value for each position in the domain of $f$ using 61 samples per location. This is rather small for the Monte Carlo method but still produces a nice approximation. The norm of the difference between the two numerical adjoint solutions is $0.0443$.

Finally, we reused the trajectories produced from this adjoint simulation to simulate the Boltzmann quantity via \eqref{eq:boltz_sde_solution}. The results are shown in Fig.~\ref{fig:boltz_reuse}. The norm of the difference between Fig.~\ref{fig:boltz_reuse} and the traditional Boltzmann simulation in Fig.~\ref{fig:boltz_trad} was $0.0162$.

\subsection{Production-Level Example}
\label{sec:its}

\begin{figure}[htp]
\centering
\begin{subfigure}{.49\textwidth}
\includegraphics[width=\textwidth]{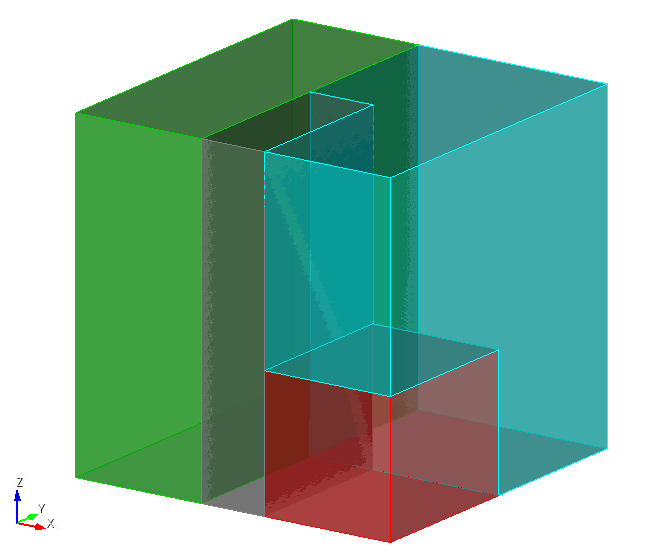}
\caption{Problem Geometry}
\label{fig:geo}
\end{subfigure}
\begin{subfigure}{.49\textwidth}
\includegraphics[width=\textwidth]{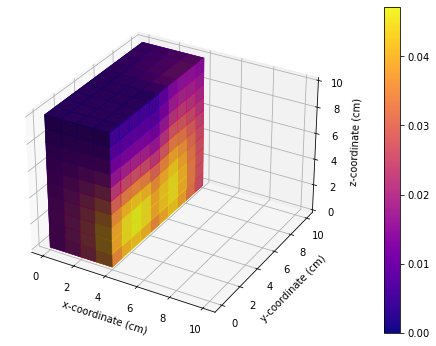}
\caption{Boltzmann Fluence}
\label{fig:fluence_its}
\end{subfigure}
\begin{subfigure}{.49\textwidth}
\includegraphics[width=\textwidth]{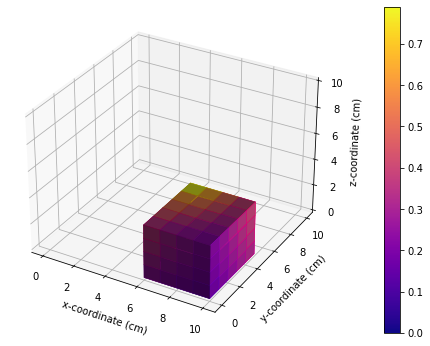}
\caption{Adjoint, Reused Trajectories}
\label{fig:its_adj_reuse}
\end{subfigure}
\begin{subfigure}{.49\textwidth}
\includegraphics[width=\textwidth]{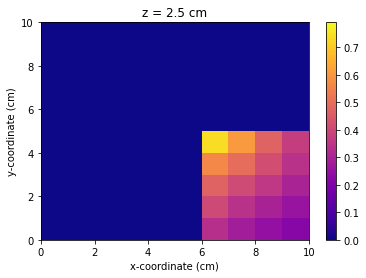}
\caption{Adjoint, $z=2.5$.}
\label{fig:its_adj_slice}
\end{subfigure}
\caption{\textbf{(a)} Geometry of ITS problem. Red region represents a volumetric source of gamma radiation.  The detector region is illustrated in green, a lead wall is illustrated in gray, and the remaining blue region represents water. Both the source and detector regions are also taken to be water.  \textbf{(b)} Boltzmann fluence simulated in ITS. This is calculated by summing the individual detector response from all simulated point sources. \textbf{(c)} Adjoint quantity calculated from reusing trajectories. \textbf{(d)} Single slice of \textbf{(c)} about $z=2.5$.}
\label{fig:numerical_ex2}
\end{figure}

To further demonstrate the practicality of the method, we use a simulation in the ITS production code.  Here we simulate a photon transport problem with continuous-energy cross sections.  A volumetric source of gamma radiation is located in the corner of a block of material.  The source region is illustrated in red in Figure~\ref{fig:geo}. The source is distributed uniformly in the $4 \times 5 \times 4$ cm region. Emitted photons are sampled uniform in energy from 1-5 MeV and uniformly isotropic in angle. A detector region, illustrated in green, is $4 \times 10 \time 10$ cm. The detector is partially separated from the source by a $2$ cm thick lead wall illustrated in gray. The source, detector, and remainder of the block, illustrated in blue, are all modeled as water. The block is surrounded by vacuum boundary conditions.  For simplicity, photons were only simulated in the energy range of 1-5 MeV.

For convenience, the simulation was carried out over a series of point sources representing a discretization of the source by starting locations.  The fluence is reported for each of the locations in the detector region. Each starting location is initialized with $1\times10^8$ particles, with uniform selections in energy and uniformly isotropic selections in angle. Cumulative particle counts for each of these detector regions can be backed out of the fluence values by multiplying by the source strength, the volume of the voxel --- 1 by construction, and the total number of particles. These counts can be treated as the results of a particle passing through an indicator function $g$ as in \eqref{eq:adj_solution}. The absorption term is handled in the simulation of the stochastic particles by removing them from simulation once they are absorbed.  Hence, we need only average the appropriate sums of these tallies to obtain the adjoint quantity for each of the starting locations represented by the point source array.  Note, the information we obtained from the ITS solution was integrated in both angle and energy; therefore the adjoint quantity we obtain will also be integrated in angle and energy.

We plot our results of this reuse of trajectory information in Figure \ref{fig:its_adj_reuse}.  Furthermore, we plot a single slice of the simulation about $z=2.5$cm in Figure \ref{fig:its_adj_slice}.  This would be the $3^\text{rd}$ level from the bottom when viewing Fig. \ref{fig:its_adj_reuse}.

\section{Discussion and Conclusion}
\label{sec:discussion}

Our use of SDEs has an immediate practical benefit for transport simulation because source iteration particle trajectories can be reused for adjoint transport of the expectation in \eqref{eq:adj_solution}.  Thus a single set of particle trajectories can be used to obtain both the Boltzmann flux density and the adjoint quantity, effectively doubling the efficiency of any code.  Moreover, a single trajectory set can be reused to provide the adjoint quantity for \emph{as many detectors/adjoint source terms} that are of  interest. In traditional implementations, the adjoint quantity must be simulated anew for each desired detector/adjoint source term.

These same conclusions can be drawn when interpreting interpreting the adjoint particle as an SDE.  In such a case, if one is already simulating the adjoint particle, one can reuse those adjoint particle trajectories to obtain the Boltzmann flux density.  Furthermore, those particle sets can be reused over and over to obtain the Boltzmann flux density for multiple source terms.

In our simple numerical example, Figure \ref{fig:numerical_ex}, we opted to first run the traditional Boltzmann transport equation and then retain those trajectories. Then, we iterate through those data to generate the adjoint quantity.  However, an appropriately organized code could calculate all needed tallies simultaneously, allowing one to simulate Boltzmann and adjoint quantities simultaneously on the same set of trajectories.

The ITS code has an adjoint calculation mode based on multigroup cross sections. The techniques demonstrated here open a path to obtaining adjoint information from the more accurate forward continuous-energy simulations. Further, one is able to obtain both forward and adjoint information from a single simulation. The usefulness of such information has been demonstrated in the transport community for a wide variety of techniques, such as biasing generation \cite{booth2006} and sensitivity analysis \cite{perfetti2012}. This work suggests an approach of using stochastic calculus to derive such methods.

To conclude we note that beyond simulation, a stochastic calculus representation of Boltzmann transport opens up analytic tools for simulation previously unavailable.  It is conceivable that sensitivities to parameter values within cross sections or distributions could be calculated through a single sample set with the appropriate analysis.  The Malliavin calculus provides derivatives which may be of use here and are under active investigation \cite{bond2022sensitivity}.  In our immediate future work, we are actively seeking stochastic calculus representations and theorems for fission and the branching process.  A successful representation will provide further utility for these formal stochastic methods within the particle transport field.


\pagebreak
\section*{Acknowledgments}

This article has been authored by an employee of National Technology \& Engineering Solutions of Sandia, LLC under Contract No. DE-NA0003525 with the U.S. Department of Energy (DOE). The employee owns all right, title and interest in and to the article and is solely responsible for its contents. The United States Government retains and the publisher, by accepting the article for publication, acknowledges that the United States Government retains a non-exclusive, paid-up, irrevocable, world-wide license to publish or reproduce the published form of this article or allow others to do so, for United States Government purposes. The DOE will provide public access to these results of federally sponsored research in accordance with the DOE Public Access Plan \url{https://www.energy.gov/downloads/doe-public-access-plan}. This paper describes objective technical results and analysis. Any subjective views or opinions that might be expressed in the paper do not necessarily represent the views of the U.S. Department of Energy or the United States Government.

\pagebreak

\section*{Appendix}

In this section, we provide proofs for Theorem \ref{thm:adj} and Corollary \ref{thm:boltz}. These proofs require appealing strongly to the stochastic calculus. The proofs are constructed with helpful remarks for the ease of the uninitiated reader. For an introduction to SDEs, see \cite{wiersema2008brownian}; for an in-depth treatment of stochastic calculus and integration, see \cite{protter1992stochastic}; for generalized stochastic calculus with jump diffusions, see \cite{hanson2007applied}.

Before proving our theorems, we make a few remarks.  First, we address the continuity of the SDE process at the boundary. Since the probability of a scattering event occurring exactly at the boundary is zero, the process $X(t)$ is continuous at the boundary almost surely. However, in any discretized simulation setting, it might be possible that a scattering event is generated when a particle hits the boundary. If this does occur in simulation, it is imperative to assume that particle hits the boundary first, ensuring the process is continuous at the boundary point.

Next, these proofs are given in the single spatial dimension case.  However, the proofs for all three spatial dimensions and two direction dimensions is identical, requiring only three additional SDEs in \eqref{eq:adj_sde} and the additional terms arising from It\={o}'s rule in \eqref{eq:ito_aux}.  For time-dependent problems with either an initial or final condition, the proof follows as the Feynman-Kac formula for jump-diffusions (see \cite{hanson2007applied}, Ch. 7).

Finally, we note that our results do not prove that solutions to the equations exist.  Rather, we only include the assumptions necessary to write the solutions in probabilistic form. Existence of solutions is well studied and often involve balancing conditions with the boundary condition; one sufficient condition for a solution to the Boltzmann equation with any boundary condition is that the cross sections are bounded and the source term is square-integrable \cite{pao1973nonlinear}. These are similar to our conditions. We will also remark that our assumptions were geared toward a dominated convergence argument and that it may be possible to weaken our assumptions so that a similar limit in $n$ could be taken with a uniformly integrable argument.

\begin{proof}[Proof (of Theorem \ref{thm:adj})]
We require the use of the Poisson random measure $\mathcal{N}$. The measure can be related to the counting process $\rmd N$ as:
\begin{equation*}
\int_{\eta,\theta\in\mathcal{A}}\mathcal{N}\left(\left(t+\rmd t\right],\left(\eta + \rmd\eta\right],\left(\theta + \rmd\theta\right],X(t),t\right) = \rmd N\left(t,\omega,\varepsilon,X(t)\right).
\end{equation*}
Here, it is assumed that the random variables $\omega$ and $\varepsilon$ (and indeed the deterministic coordinate variables $\Omega$ and $E$) live in some space $\mathcal{A} = [-1,1]\times\mathcal{E}$ and the integration is performed in this space; the value $\eta$ can be any realization of the random variable $\omega$ and the value $\theta$ can be any realization of the random variable $\varepsilon$. The Poisson random measure can be written in terms of a mean-zero Poisson random measure, $\widehat{\mathcal{N}}$, and its mean:
\begin{align}
\begin{split}
\mathcal{N}&\left(\left(t+\rmd t\right], \left(\eta +\rmd\eta\right],\left(\theta + \rmd \theta\right],X(t),t\right)\\
&= \widehat{\mathcal{N}}\left(\left(t+\rmd t\right], \left(\eta +\rmd\eta\right],\left(\theta + \rmd \theta\right],X(t),t\right) + \Sigma_s\left(X(t)\right)p\left(\eta,\theta\,\middle|\,X(t)\right)\,\rmd\eta\rmd\theta\rmd t.
\end{split}
\label{eq:zero_measure}
\end{align}

This proof requires the use of an auxiliary function.  Consider the function
\begin{equation*}
w\left(X(t)\right) = \Psi\left(X(t)\right)\exp\left(-\int_0^t\Sigma_a\left(X(u)\right)\rmd u\right),
\end{equation*}
and for $n\in\mathbb{N}$ define
\begin{equation*}
T_{x,\Omega,E,n} = \min\left\{n,\tau_{x,\Omega,E}\right\}.
\end{equation*}
For compactness, let $H\left(X(t),\omega,\varepsilon\right)=(X_1(t),\omega,\varepsilon)^\top-X(t)$. Given an $n$, for $t\leq T_{x,\Omega,E,n}$ It\={o}'s rule (a chain rule for stochastic calculus) yields
\begin{align}
\begin{split}
\rmd &w\left(X(t)\right)=\exp\left(-\int_0^t\Sigma_a\left(X(u)\right)\rmd u\right)\left[\left(X_2(t)\frac{\partial}{\partial x}\Psi\left(X(t)\right) - \Sigma_a\left(X(t)\right)\Psi\left(X(t)\right)\right)\rmd t\right.\\
&+\left.\int_{\eta,\theta\in\mathcal{A}}\left(\Psi\left(X(t)+H\left(X(t),\eta,\theta\right)\right)-\Psi\left(X(t)\right)\right)\mathcal{N}\left(\left(t+\rmd t\right],\left(\eta+\rmd \eta\right],\left(\theta+\rmd \theta\right],X(t),t\right)\right],\\
&=\exp\left(-\int_0^t\Sigma_a\left(X(u)\right)\rmd u\right)\left[-g\left(X(t)\right)\rmd t\phantom{\int}\right.\\
&+\left.\int_{\eta,\theta\in\mathcal{A}}\left(\Psi\left(X_1(t),\eta,\theta\right)-\Psi\left(X(t)\right)\right)\widehat{\mathcal{N}}\left(\left(t+\rmd t\right],\left(\eta+\rmd \eta\right],\left(\theta+\rmd\theta\right],X(t),t\right)\right].
\end{split}
\label{eq:ito_aux}
\end{align}
The second equality is obtained by expanding $\mathcal{N}$ via \eqref{eq:zero_measure} and then by employing \eqref{eq:adj_bvp}. Now, we integrate both sides of this result from $0$ to $T_{x,\Omega,E,n}$ and take an expectation conditioning on $X(0)=\left(x,\Omega,E\right)^\top$. Since the expected value of an integral against a mean-zero Poisson random measure is zero, this leaves us with
\begin{align*}
\mathbb{E}&\left[w\left(X\left(T_{x,\Omega,E,n}\right)\right)\,\middle|\,X(0)=\left(x,\Omega,E\right)^\top\right] - \mathbb{E}\left[w\left(X(0)\right)\,\middle|\,X(0)=\left(x,\Omega,E\right)^\top\right]\\
&=\mathbb{E}\left[-\int_0^{T_{x,\Omega,E,n}}g\left(X(t)\right)\exp\left(-\int_0^t\Sigma_a\left(X(u)\right)\rmd u\right)\rmd t\,\middle|\,X(0)=\left(x,\Omega,E\right)^\top\right].
\end{align*}
The second expectation on the left-hand side simplifies to $w\left(x,\Omega,E\right)=\Psi\left(x,\Omega,E\right)$. Solving for $\Psi$ yields
\begin{align}
\begin{split}
\Psi\left(x,\Omega,E\right) &=\mathbb{E}\left[w\left(X\left(T_{x,\Omega,E,n}\right)\right),\middle|\,X(0)=\left(x,\Omega,E\right)^\top\right]\\
&+\mathbb{E}\left[\int_0^{T_{x,\Omega,E,n}}g\left(X(t)\right)\exp\left(-\int_0^t\Sigma_a\left(X(u)\right)\rmd u\right)\rmd t\,\middle|\,X(0)=\left(x,\Omega,E\right)^\top\right].
\end{split}
\label{eq:psi_almost}
\end{align}
We would like to take the limit as $n\to\infty$ of this result. Since $\mathbb{E}\left[\tau_{x,\Omega,E}\right]<\infty$, as $n\to\infty$, we will have $T_{x,\Omega,E,n}\to\tau_{x,\Omega,E}$. Since $X$ is continuous at the boundary, $X\left(T_{x,\Omega,E,n}\right)\to X\left(\tau_{x,\Omega,E}\right)$. In order to take the desired limit and push the limit inside both expectations, by dominated convergence it is sufficient to show that the interiors of each expectation are bounded for all $n$.

Beginning with the first term, note that for a given $t<T_{x,\Omega,E,n}$ that $X_2$ and $X_3$ (and hence $X$) must have finitely many discontinuities.  If they did not, then there were infinite arrivals of a Poisson process with a bounded rate $\Sigma_s$ in a finite time period, a contradiction. Next, since $\Sigma_a$ is continuous almost everywhere and since $X$ has only finitely many discontinuities, $\Sigma_a(X(t))$ must also be continuous almost everywhere and hence integrable. Therefore $\exp\left(-\int_0^{T_{x,\Omega,E,n}}\Sigma_a\left(X(t)\right)\rmd t\right)$ is continuous. Furthermore, this term is bounded since $\Sigma_a$ is bounded and since $T_{x,\Omega,E,n}$ is finite for all $n$. Since $\Psi$ is a classical solution to \eqref{eq:adj_bvp}, it is also continuous. Since $\mathcal{D}$ is compact, $\Psi$ is also bounded. Hence $w$ is continuous and $w\left(X\left(T_{x,\Omega,E,n}\right)\right)$ is bounded for all $n$. Ergo:
\begin{align*}
\lim_{n\to\infty}\mathbb{E}&\left[w\left(X\left(T_{x,\Omega,E,n}\right)\right)\,\middle|\,X(0)=\left(x,\Omega,E\right)^\top\right],\\
&=\mathbb{E}\left[w\left(X\left(\tau_{x,\Omega,E}\right)\right)\,\middle|\,X(0)=\left(x,\Omega,E\right)^\top\right],\\
&=\left[B\left(X\left(\tau_{x,\Omega,E}\right)\right)\exp\left(-\int_0^{\tau_{x,\Omega,E}}\Sigma_a\left(X(t)\right)\rmd t\right)\,\middle|\,X(0)=\left(x,\Omega,E\right)^\top\right].
\end{align*}

Next we consider the integral term.  For $t\leq T_{x,\Omega,E,n}$, similar arguments to the proceeding show that $\exp\left(-\int_0^t\Sigma_a\left(X(u)\right)\rmd u\right)$ is both continuous and bounded. Since $g$ is continuous almost everywhere, their product is also continuous almost everywhere and therefore integrable. Since $g$ is also bounded and since $T_{x,\Omega,E,n}<\infty$ for all $n$, the integral in the expectation in \eqref{eq:psi_almost} is both continuous in $n$ and bounded for all $n$. Hence:
\begin{align*}
\lim_{n\to\infty}\mathbb{E}&\left[\int_0^{T_{x,\Omega,E,n}}g\left(X(t)\right)\exp\left(-\int_0^t\Sigma_a\left(X(u)\right)\rmd u\right)\rmd t\,\middle|\,X(0)=\left(x,\Omega,E\right)^\top\right]\\
&=\mathbb{E}\left[\int_0^{\tau_{x,\Omega,E,}}g\left(X(t)\right)\exp\left(-\int_0^t\Sigma_a\left(X(u)\right)\rmd u\right)\rmd t\,\middle|\,X(0)=\left(x,\Omega,E\right)^\top\right].
\end{align*}
Therefore, taking the limit in $n$ of \eqref{eq:psi_almost} produces \eqref{eq:adj_sol_general}, proving the result.
\end{proof}

\begin{proof}[Proof (of Corollary \ref{thm:boltz})]
The proof rests on rearranging \eqref{eq:boltz_problem} to fit the structure of \eqref{eq:adj_bvp}.  Up to constant multiples, the only term out of place is the integration in the probability distribution. Using \eqref{eq:q_def}, we have \small
\begin{equation*}
\int\Phi\left(x,\Omega',E'\right)\Sigma_s\left(x,\Omega',E'\right)p\left(\Omega, E\,\middle|\,\Omega',E'\right)\rmd E'\rmd\Omega' = \int\Phi\left(x,\Omega',E'\right)S\left(x,\Omega,E\right)q\left(\Omega',E'\,\middle|\,\Omega, E\right)\rmd E'\rmd\Omega'.
\end{equation*}\normalsize
Since $\Sigma_s$ is bounded and continuous almost everywhere, so is $S$ by \eqref{eq:back_scatter}. Similarly, $\Sigma_a+\Sigma_s-S$ is bounded and continuous almost everywhere. Hence Theorem \ref{thm:adj} applies and the result follows.
\end{proof}

\end{document}